\documentclass[12pt]{article}
\usepackage{amsmath,amssymb,euscript,amsthm,amsfonts}
\usepackage[cp1251]{inputenc}
\usepackage[english,russian,ukrainian]{babel}

\usepackage{color}
\definecolor{darkblue}{rgb}{0.00,0.25,0.50}
\usepackage[colorlinks,filecolor=blue,citecolor=darkblue]{hyperref}
\usepackage{cite}
\begin{document}

\begin{center}
\textbf{\Large Оцінки рівномірних наближень сумами Зигмунда на
класах згорток періодичних функцій}
\end{center}
\vskip0.5cm

\begin{center}
А.\,С.~Сердюк \\ \emph{\small Інститут математики НАН України,
Київ} \\
У.\,З.~Грабова \\ \emph{\small Східноєвропейський національний
університет імені Лесі Українки, Луцьк}
\end{center}
\vskip0.5cm

\begin{abstract} We obtain order-exact estimates for uniform
approximations by using Zygmund sums $Z^{s}_{n}$ of classes
$C^{\psi}_{\beta,p}$ of $2\pi$-periodic continuous functions $f$
representable by convolutions  of functions from unit balls of the
space  $L_{p}$, $1< p<\infty$, with a fixed kernels
$\Psi_{\beta}\in L_{p'}$, $\frac{1}{p}+\frac{1}{p'}=1$. In
addition, we  find a set of allowed values of parameters (that
define the class $C^{\psi}_{\beta,p}$ and the linear method
$Z^{s}_{n}$) for which Zygmund sums and Fejer sums realize the
order of the best uniform approximations by trigonometric
polynomials of those classes.

\vskip0.6cm

Одержано точні за порядком оцінки рівномірних наближень сумами
Зигмунда $Z^{s}_{n}$ на класах $C^{\psi}_{\beta,p}$
 $2\pi$-періодичних неперервних функцій $f$, які зображуються у вигляді згортки функцій, що належать одиничним кулям просторів $L_{p}$, $1< p<\infty$,
  з фіксованими твірними ядрами $\Psi_{\beta}\in  L_{p'}$,
  $\frac{1}{p}+\frac{1}{p'}=1$.
Вказано множину допустимих значень параметрів (що визначають класи
$C^{\psi}_{\beta,p}$ та лінійний метод $Z^{s}_{n}$) при яких суми
Зигмунда, а також суми Фейєра, забезпечують порядок найкращих
рівномірних наближень тригонометричними поліномами на вказаних
класах.
\end{abstract}

\vskip1cm

 \normalsize \vskip 3mm
Нехай $L_{p}$, $1\leq p\leq\infty$ --- простір $2\pi$-періодичних
сумовних функцій $f$ зі скінченною нормою $\| f\|_{p}$, де при
$p\in [1,\infty)$
\begin{equation*}
\|
f\|_{p}=\Big(\int\limits_{0}^{2\pi}|f(t)|^{p}dt\Big)^{\frac{1}{p}},
\end{equation*}
а при $p=\infty$ \
\begin{equation*}\|f\|_{\infty}=\mathop{\rm{ess}\sup}\limits_{t}|f(t)|;
\end{equation*} $C$
--- простір $2\pi$-періодичних неперервних функцій, у
якому норма задається рівністю $\|f\|_{C}=\max\limits_{t}|f(t)|$.

Нехай, далі,  $L^{\psi}_{\beta,p}$, \ $1\leq p<\infty$ --- клас
$2\pi$-періодичних функцій $f(x)$, що майже для всіх
$x\in\mathbb{R}$ представляються згортками

\hskip -5mm\begin{equation}\label{zgo}
f(x):=\frac{a_{0}}{2}+(\Psi_{\beta}\ast\varphi)(x)=\frac{a_{0}}{2}+
\frac{1}{\pi}\int\limits_{-\pi}^{\pi}\Psi_{\beta}(x-t)\varphi(t)dt,
\ a_{0}\in\mathbb{R}, \ \varphi\perp 1,
\end{equation}
де  $\|\varphi \|_{p}\leq 1$,  а $\Psi_{\beta}(t)$ --- сумовна на
$[0,2\pi)$ функція, ряд Фур'є якої має вигляд
\begin{equation}\label{rff}
\sum\limits_{k=1}^{\infty}\psi(k)\cos
\big(kt-\frac{\beta\pi}{2}\big), \ \psi(k)>0,\ \beta\in
    \mathbb{R}.
\end{equation}
Функцію $\varphi$ в зображенні \eqref{zgo}, згідно з О.І.
Степанцем \cite{S1}, с. 132), називають $(\psi,\beta)$-похідною
функції $f$ і позначають через $f^{\psi}_{\beta}$.

 Якщо $\psi(k)=k^{-r}$, $r>0$, $\beta\in  \mathbb{R}$, то класи $L^{\psi}_{\beta,p}$ є  класами
Вейля-Надя і  позначаються через $W^{r}_{\beta,p}$.  При
$\beta=r$, \ $r\in \mathbb{N}$ останні класи є відомими класами
$W^{r}_{p}$.

Якщо твірне ядро $\Psi_{\beta}$ класу $L^{\psi}_{\beta,p}$, $1\leq
p\leq\infty$ задовольняє включенню $\Psi_{\beta}\in L_{p'}$,
\linebreak$\frac{1}{p}+\frac{1}{p'}=~1$, то
$L^{\psi}_{\beta,p}\subset L_{\infty}$,  а згортки виду
\eqref{zgo} є неперервними функціями (див. твердження 3.8.1 та
3.8.2 роботи \cite{S1}, с. 137, 138). Тому клас усіх функцій $f$
виду \eqref{zgo}, для яких $\|\varphi\|_{p}\leq 1$,
$\Psi_{\beta}\in L_{p'}$  будемо позначати через
$C^{\psi}_{\beta,p}$. Зрозуміло, що у випадку $f\in
C^{\psi}_{\beta,p}$ рівність \eqref{zgo} виконується при всіх
$x\in\mathbb{R}$.

Вважаючи, що послідовність $\psi(k)$, що визначає клас
$C^{\psi}_{\beta,p}$, є слідом на множині $\mathbb{N }$ деякої
неперервної функції $\psi(t)$ неперервного аргументу $t\geq1$,
позначимо через $\Theta_{p}$, \linebreak$1\leq p<\infty$, множину
монотонно незростаючих функцій $\psi(t)$, для яких існує стала
\linebreak$\alpha>\frac{1}{p}$ така, що функція
$t^{\alpha}\psi(t)$ майже спадає, тобто знайдеться додатна стала
$K$ така, що $t^{\alpha}_{1}\psi(t_{1})\leq
Kt^{\alpha}_{2}\psi(t_{2})$ для будь-яких $t_{1}>t_{2}\geq 1$.
 Умова $\psi\in \Theta_{p}$, як неважко переконатись,
гарантує справедливість включення $\Psi_{\beta}\in L_{p'}$,
$\frac{1}{p}+\frac{1}{p'}=1$ (див., наприклад \cite{Bari}, с.
657).

Прикладами функцій $\psi$, що задовольняють умову  $\psi\in
\Theta_{p}$ є, зокрема, функції виду
$\psi_{1}(t)=\frac{1}{t^{r}}$, \ $r>\frac{1}{p}$;
$\psi_{2}(t)=\frac{\ln^{\alpha}(t+c)}{t^{r}}$, \ $r>\frac{1}{p}$,
\ $\alpha>0$, \ $c>e^{\frac{2\alpha}{r-\frac{1}{p}}}-1$;
$\psi_{3}(t)=\frac{1}{t^{r}\ln^{\alpha}(t+c)}$, \ $r>\frac{1}{p}$,
\ $\alpha>0$, \ $c>0$; \  $\psi_{4}(t)=\frac{\ln
\ln^{\alpha}(t+c)}{t^{r}}$, \ $r>\frac{1}{p}$, \ $\alpha>0$, \
$c>e^{\frac{2\alpha}{r-\frac{1}{p}}}-1$.

 Сумами Зигмунда
 функції $f\in L_{1}$ називають тригонометричні
поліноми виду
\begin{equation}\label{sz}
Z_{n}^{s}(f;t)=\frac{a_{0}}{2}+\sum\limits_{k=1}^{n-1}\left(1-\left(\frac{k}{n}\right)^{s}\right)(a_{k}(f)\cos
kt+b_{k}(f)\sin kt), s>0,
\end{equation}
де $a_{k}(f)$ і $b_{k}(f)$ --- коефіцієнти Фур'є функції $f$. Суми
Зигмунда при довільному $s>0$ були введені А. Зигмундом
\cite{Zi1}.

При $s=1$ суми $Z^{s}_{n}$ перетворюються у відомі суми Фейєра
\begin{equation}\label{sf}
\sigma_{n}(f;t)=\frac{a_{0}}{2}+\sum\limits_{k=1}^{n-1}\left(1-\frac{k}{n}\right)(a_{k}(f)\cos
kt+b_{k}(f)\sin kt),
\end{equation}
де $a_{k}(f)$ і $b_{k}(f)$ --- коефіцієнти Фур'є функції $f$.

При всіх $1\leq p$, $q\leq\infty$ і $\beta\in
    \mathbb{R}$ таких, що $L^{\psi}_{\beta,p}\subset L_{q}$
    розглянемо величини вигляду
\begin{equation}\label{velzag}
{\cal E}_{n}\left(L^{\psi}_{\beta,p};
Z_{n}^{s}\right)_{{q}}=\mathop{\sup}\limits_{f\in
L^{\psi}_{\beta,p}}\|f(\cdot)-Z^{s}_{n}(f;\cdot)\|_{{q}}.
\end{equation}

 В роботі досліджуються порядкові оцінки  величин \eqref{velzag}
 при $\beta\in
    \mathbb{R}$, $\psi\in \Theta_{p}$, $1< p<\infty$ і
    $q=\infty$. Зрозуміло, що при $\psi\in \Theta_{p}$, $1<
    p<\infty$,
    $q=\infty$
\begin{equation}\label{pv1}
{\cal E}_{n}\left(L^{\psi}_{\beta,p};
Z_{n}^{s}\right)_{\infty}={\cal E}_{n}\left(C^{\psi}_{\beta,p};
Z_{n}^{s}\right)_{C}=\mathop{\sup}\limits_{f\in
C^{\psi}_{\beta,p}}\|f(\cdot)-Z^{s}_{n}(f;\cdot)\|_{C}.
\end{equation}

Порядкові оцінки величин ${\cal
E}_{n}\left(W^{r}_{\infty};Z_{n}^{s}\right)_{C}$ для  $r\in
    \mathbb{N}$  знайдено в  роботі А. Зигмунда \cite{Zi1}.
 Б. Надь \cite{N} дослідив
 величину ${\cal
E}_{n}\left(W^{r}_{\beta,\infty};Z_{n}^{s}\right)_{C}$ при $r>0$,
$\beta\in \mathbb{Z}$,  причому для $s\leq r$ ним знайдено для цих
величин асимптотичні рівності, а для $s> r$ --- порядкові оцінки.
Згодом С.А.~ Теляковський \cite{Telyakovskiy63} для $r>0$ одержав
асимптотично точні   рівності для  величин ${\cal
E}_{n}\left(W^{r}_{\beta,\infty};Z_{n}^{s}\right)_{C}$ при
$\beta\in \mathbb{R}$, $n\rightarrow\infty$.

Для сум Фейєра $\sigma_{n}(f;t)$ порядкові оцінки величин ${\cal
E}_{n}\left(W^{r}_{\beta,p};\sigma_{n}\right)_{q}$ при $\beta\in
\mathbb{Z}$ і $p=q=\infty$ були знайдені С.М. Нікольським
\cite{Nik}; для величин ${\cal
E}_{n}\left(W^{r}_{p};\sigma_{n}\right)_{q}$ у випадках а)$1<p, \
q<\infty$; \  \  б)$1\leq p<\infty, \ q=\infty$; \ \ в)$p=1, \
1<q<\infty$ ---
 у роботах В.М. Тихомирова \cite{TUX} та А.І.~ Камзолова
\cite{Kam1}, \cite{Kam2}. Згодом результати \cite{TUX} ---
\cite{Kam2} були доповнені у роботах М.В. Костича \cite{Kostuch1},
\cite{Kostuch2}, де  знайдено порядкові оцінки величин ${\cal
E}_{n}\left(W^{r}_{\beta,p};Z_{n}^{s}\right)_{q}$, $\beta\in
\mathbb{R}$ при
 $1\leq p<\infty$, $q=\infty$, $r>\frac{1}{p}$, а також при $p=1$,
    $1<q<\infty$, $r>1-\frac{1}{p}$.

На класах $C^{\psi}_{\beta,p}$ та $L^{\psi}_{\beta,p}$ оцінки
наближень сумами Зигмунда в рівномірній та інтегральних метриках
вивчались у роботах  Д.М. Бушева \cite{Bu}, І.Б. Ковальської
\cite{Koval} та інших. В роботі \cite{Koval} у випадку  $1<q\leq
p<\infty$ одержано, зокрема, точні порядкові оцінки величин ${\cal
E}_{n}\left(L^{\psi}_{\beta,p};Z_{n}^{s} \right)_{q}$  за умови
монотонного спадання до нуля послідовності $\psi(k)$, а у випадку
$1<p<q<\infty$ --- за умови $\frac{\psi(k)}{\psi(ck)}\leq
K<\infty$, $c>1$, $K>0$  і монотонного спадання до нуля
послідовності $\psi(k)k^{\frac{1}{p}-\frac{1}{q}}$. Якщо $p=q=2$ в
\cite{Koval} знайдено точні рівності величин вигляду
\eqref{velzag}. Точні рівності для вказаних величин при $p=2$ і
$q=\infty$, а також при $p=1$ і $q=2$ за умови збіжності ряду
$\sum\limits_{k=1}^{\infty}\psi^{2}(k)$ встановлено у роботах А.С.
Сердюка та І.В. Соколенка \cite{Serduk1}, \cite{Serduk2}.

Асимптотично точні оцінки величин ${\cal
E}_{n}\left(L^{\psi}_{\beta,p};Z^{s}_{n}\right)_{q}$ при деяких
досить природних обмеженнях на функції $\psi$ знайдено у роботах
Д.М.~ Бушева \cite{Bu} при $p=q=\infty$ та І.Б.~ Ковальської
\cite{Koval} при $p=q=1$.

Поряд з величинами \eqref{pv1} в роботі  розглядаються також
найкращі рівномірні наближення класів $C^{\psi}_{\beta,p}$
тригонометричними поліномами  $t_{n-1}$ порядку $n-1$, тобто
величини вигляду
\begin{equation}\label{nnabl}
{E}_{n}(C^{\psi}_{\beta,p})_{C}=\sup\limits_{f\in
C^{\psi}_{\beta,p}}\inf\limits_{t_{n-1}}\|f(\cdot)-t_{n-1}(\cdot)\|_{C},
\ \ 1< p<\infty.
\end{equation}
Дослідженню порядкових рівностей для величин \eqref{nnabl}
присвячено роботу \cite{Gra} (там же можна більш детально
ознайомитись із наявною бібліографією щодо оцінок вказаних
величин).

В даній   роботі  знайдено порядкові рівності для величин
\eqref{pv1} при довільних  $1< p<\infty$, \ $q=\infty$ і $\psi\in
\Theta_{p}$ і тим самим доповнено встановлені  в \cite{Koval}
оцінки величин \eqref{velzag} на випадки $1< p<\infty$, \
$q=\infty$, а також
 узагальнено  результати роботи \cite{Kostuch1} на класи
  $C^{\psi}_{\beta,p}$, $\psi\in \Theta_{p}$,
$\beta\in\mathbb{R}$. Крім того,  в роботі встановлено множину
допустимих значень параметрів (що визначають класи
$C^{\psi}_{\beta,p}$, та лінійний метод $Z^{s}_{n}$) при яких суми
Зигмунда, а також суми Фейєра, забезпечують порядок найкращих
рівномірних наближень тригонометричними поліномами на вказаних
класах, тобто
\begin{equation*}
{\cal E}_{n}\left(C^{\psi}_{\beta,p};
Z_{n}^{s}\right)_{C}\asymp{E}_{n}(C^{\psi}_{\beta,p})_{C}.
\end{equation*}
 Перейдемо до точних
формулювань.

Будемо казати, що додатна функція $g(t)$, задана на $[1,\infty)$,
належить до множини $A^{+}$ ( і записувати $g\in A^{+}$), якщо
існує $\varepsilon>0$ таке, що $g(t)t^{-\varepsilon}$ зростає на
$[1,\infty)$. Аналогічно,  якщо існує $\varepsilon>0$ таке, що
$g(t)t^{\varepsilon}$ спадає на $[1,\infty)$, то будемо казати, що
$g$ належить до множини $A^{-}$ і записувати $g\in A^{-}$.

Через $\mathcal{Z}$ позначимо множину неперервних слабо коливних
(в сенсі Зигмунда) функцій, тобто  додатних функцій $g(t)$,
визначених на $[1,\infty)$, таких, що при довільному $\delta>0$
функція $g(t)t^{\delta}$ зростає, а $g(t)t^{-\delta}$ спадає, для
достатньо великих $t$. При фіксованому $\rho>1$ через
$\mathcal{Z}^{+}_{\rho}$ позначимо підмножину монотонно зростаючих
функцій $g\in \mathcal{Z}$, що задовольняють умову
\begin{equation*}
 g^{\rho}(n)\ln
n=O\Big(\int\limits_{1}^{n}\frac{g^{\rho}(t)}{t}dt\Big), \ n\in
\mathbb{N},
\end{equation*}
а через $\mathcal{Z}^{-}_{\rho}$ --- підмножину монотонно спадних
функцій $g\in \mathcal{Z}$, що задовольняють умову
\begin{equation*}
\int\limits_{1}^{n}\frac{g^{\rho}(t)}{t}dt=O\big(g^{\rho}(n)\ln
n\big), \ n\in \mathbb{N}.
\end{equation*}

 Під записом $A(n)=O(B(n))$ розуміємо, що
існує стала $K>0$, така, що виконується нерівність $A(n)\leq
K(B(n))$, де $n\in \mathbb{N}$. Запис $A(n)\asymp B(n)$ означає,
що $A(n)=O(B(n))$ і одночасно $B(n)=~O(A(n))$.

\vskip 3mm

\bf Теорема 1.  \it{Нехай $1< p<\infty$, $s>0$ і
$g_{s,p}(t):=\psi(t)t^{s+\frac{1}{p}}$. Тоді

1. Якщо $\psi\in \Theta_{p}$ і $g_{s,p}\in A^{+}$,  то  для
довільних $\beta\in \mathbb{R}$ і $n\in \mathbb{N}$
\begin{equation}\label{t5}
{\cal E}_{n}\left(C^{\psi}_{\beta,p};
Z_{n}^{s}\right)_{C}=O\big(\psi(n)n^{\frac{1}{p}}\big).
\end{equation}

 2. Якщо   $g_{s,p}\in \mathcal{Z}$, то  для довільних
$\beta\in \mathbb{R}$ і $n\in \mathbb{N}$

\begin{equation}\label{t4}
{\cal E}_{n}\left(C^{\psi}_{\beta,p};
Z_{n}^{s}\right)_{C}=O\Bigg(\frac{1}{n^{s}}\Bigg(\int\limits_{1}^{n}\frac{\big(\psi(t)t^{s+\frac{1}{p}}\big)^{p'}}{t}dt\Bigg)^{\frac{1}{p'}}\Bigg),
\ \frac{1}{p}+\frac{1}{p'}=1.
\end{equation}

 3. Якщо $g_{s,p}\in A^{-}$, то  для довільних
$\beta\in \mathbb{R}$ і $n\in \mathbb{N}$
\begin{equation}\label{t1}
{\cal E}_{n}\left(C^{\psi}_{\beta,p};
Z_{n}^{s}\right)_{C}=O\Big(\frac{1}{n^{s}}\Big).
\end{equation}

{\textbf{\textit{Доведення. }} \ \rm Оскільки оператор
$$Z^{s}_{n}: f(t)\rightarrow Z^{s}_{n}(f,t)$$
 є лінійним поліноміальним
оператором, інваріантним відносно зсувів
\begin{equation*}
Z^{s}_{n}(f_{h}, t)= Z^{s}_{n}(f, t+h), \ f_{h}(t)=f(t+h), \ h\in
\mathbb{R},
\end{equation*}
і норма в $C$, та класи $C^{\psi}_{\beta,p}$ також інваріантні
відносно зсувів, тобто
\begin{equation*}
\|f_{h}(t)\|_{C}=\|f(t)\|_{C}; \ f(t)\in
C^{\psi}_{\beta,p}\Rightarrow f_{h}(t)\in C^{\psi}_{\beta,p},
\end{equation*}
то
\begin{equation*}
{\cal E}_{n}\left(C^{\psi}_{\beta,p};
Z_{n}^{s}\right)_{C}=\mathop{\sup}\limits_{f\in C^{\psi}_{\beta,p}
}|f(0)-Z_{n}^{s}(f;0)|.
\end{equation*}
В силу формул \eqref{zgo} і \eqref{sz} для довільної $f\in
C^{\psi}_{\beta,p}$, $1\leq p<\infty$, $\beta\in \mathbb{R}$,
$s>0$
\begin{equation}\label{predst}
f(0)-Z_{n}^{s}(f;0)=\frac{1}{\pi}\!\!\int\limits_{-\pi}^{\pi}\!\!\!
\left(\frac{1}{n^{s}}\sum\limits_{k=1}^{n-1}\psi(k)k^{s}\cos\left(kt+\frac{\beta\pi}{2}\right)\!+\!\Psi_{-\beta,n
    }(t)\right)\!\varphi(t)dt,
\end{equation}
де $\Psi_{-\beta,n
    }(t)=\sum\limits_{k=n}^{\infty}\psi(k)\cos\left(kt+\frac{\beta\pi}{2}\right)$,
\ \  $\|\varphi\|_{p}\leq1$, $n\in \mathbb{N}$.

 З нерівності Гельдера та нерівності трикутника
випливає, що при $1\leq p<\infty$
\begin{equation}\label{111}
{\cal E}_{n}\left(C^{\psi}_{\beta,p};
Z_{n}^{s}\right)_{C}\leq\frac{1}{\pi
}\Bigg\|\frac{1}{n^{s}}\sum\limits_{k=1}^{n-1}\psi(k)k^{s}\cos\left(kt+\frac{\beta\pi}{2}\right)+
\Psi_{-\beta,n
    }(t)\Bigg\|_{p'}\leq
$$
$$
\leq\frac{1}{\pi n^{s}}\Bigg\|
\sum\limits_{k=1}^{n-1}\psi(k)k^{s}\cos\left(kt+\frac{\beta\pi}{2}\right)\Bigg\|_{p'}+\frac{1}{\pi}\big\|\Psi_{-\beta,n}(t)\big\|_{p'},
\ \frac{1}{p}+\frac{1}{p'}=1.
\end{equation}

Знайдемо оцінку зверху величини $\|\Psi_{-\beta,n}(\cdot)\|_{p'}$.
Для цього,
 застосувавши до
функції $\Psi_{-\beta,n
    }(t)$ перетворення Абеля,  одержимо
\begin{equation}\label{peretA}
\Psi_{-\beta,n
    }(t)=\sum\limits_{k=n}^{\infty}\big(\psi(k)-\psi(k+1)\big)D_{k,\beta}(t)-\psi(n)D_{n-1,\beta}(t),
\end{equation}
де
\begin{equation*}
D_{k,\beta}(t)=\frac{1}{2}\cos\frac{\beta\pi}{2}+\sum\limits_{\nu=1}^{k}\cos\left(\nu
t+\frac{\beta\pi}{2}\right)=
$$
$$
=\cos\frac{\beta\pi}{2}\Bigg(\frac{\sin\frac{2k+1}{2}t}{2\sin\frac{t}{2}}\Bigg)-\sin\frac{\beta\pi}{2}\Bigg(\frac{\cos\frac{t}{2}-\cos\frac{2k+1}
{2}t}{2\sin\frac{t}{2}}\Bigg) .
\end{equation*}
Відомо (див., наприклад, \cite{Z}, с. 89), що
\begin{equation}\label{diri}
|D_{k,\beta}(t)|=O(k), \ \
|D_{k,\beta}(t)|=O\Big(\frac{1}{|t|}\Big), \ 0<|t|\leq \pi.
\end{equation}
Тому
\begin{equation*}
\big\|
D_{k,\beta}(t)\big\|^{p'}_{p'}=\!\!\!\int\limits_{-\pi}^{\pi}\!\!|D_{k,\beta}(t)|^{p'}dt\!=O\Bigg(\int\limits_{0\leq|t|\leq\frac{1}{k}
}\!\!\!k^{p'}dt+\!\!\int\limits_{\frac{1}{k}\leq|t|\leq\pi}\!\!\!\frac{dt}{|t|^{p^{'}}}\Bigg)\!\!=\!O(k^{p'-1}).
\end{equation*}
З останньої рівності маємо при $1\leq p <\infty$, $k\in
\mathbb{N}$ і $\beta\in \mathbb{R}$
\begin{equation}\label{ndiri}
\|D_{k,\beta}(t)\|_{p'}=O(k^{1-\frac{1}{p'}})=O(k^{\frac{1}{p}}).
\end{equation}
З  \eqref{peretA} та \eqref{ndiri} одержуємо
\begin{equation}\label{npsi1}
\|\Psi_{-\beta,n}(\cdot)\|_{p'}=O\Big(\sum\limits_{k=n}^{\infty}\big(\psi(k)-\psi(k+1)\big)k^{\frac{1}{p}}\big)+\psi(n)n^{\frac{1}{p}}\Big),
\ 1\leq p<\infty.
\end{equation}
Далі, для оцінки  суми
$\sum\limits_{k=n}^{\infty}\big(\psi(k)-\psi(k+1)\big)k^{\frac{1}{p}}$
 буде корисним наступне твердження роботи \cite{Gra}.

\bf Лема 1. \it{Нехай $r\in (0,1]$, а $\psi(k)>0$, монотонно
незростає і для неї знайдеться $\varepsilon>0$ таке, що
послідовність $k^{r+\varepsilon}\psi(k)$ майже спадає.  Тоді існує
стала $K$, залежна від $\psi$ і $r$ така, що для довільних $n\in
\mathbb{N}$
\begin{equation}\label{nr2}
\psi(n)n^{r}\leq\sum\limits_{k=n}^{\infty}\big(\psi(k)-\psi(k+1)\big)k^{r}\leq
K\psi(n)n^{r}.
\end{equation}

\rm

Оскільки $\psi\in \Theta_{p}$, то, застосовуючи лему 1 при
$r=\frac{1}{p}$, із \eqref{npsi1} отримуємо оцінку
\begin{equation}\label{nb}
\|\Psi_{-\beta,n
    }(\cdot)\|_{p'}=O(\psi(n)n^{\frac{1}{p}}), \ 1\leq p<\infty.
\end{equation}

 Перейдемо до оцінки величини
\begin{equation}\label{v1}
  \Big\|
\sum\limits_{k=1}^{n-1}\psi(k)k^{s}\cos\left(kt+\frac{\beta\pi}{2}\right)\Big\|_{p'}.
\end{equation}
Застосувавши до функції
$\sum\limits_{k=1}^{n-1}\psi(k)k^{s}\cos\left(kt+\frac{\beta\pi}{2}\right)$
перетворення Абеля, одержуємо
\begin{equation*}
\sum\limits_{k=1}^{n-1}\psi(k)k^{s}\cos\Big(kt+\frac{\beta\pi}{2}\Big)=\sum\limits_{k=1}^{n-2}\Big(\psi(k)k^{s}-
\psi(k+1)(k+1)^{s}\Big)D_{k,\beta}(t)+
$$
$$
+\psi(n-1)(n-1)^{s}D_{n-1,\beta}(t)-\frac{1}{2}\cos\frac{\beta\pi}{2}.
\end{equation*}
Тоді, враховуючи \eqref{ndiri}, маємо
\begin{equation}\label{ddzn}
\Bigg\|\sum\limits_{k=1}^{n-1}\psi(k)k^{s}\cos\left(kt+\frac{\beta\pi}{2}\right)\Bigg\|_{p'}=O(1)+
$$
$$
+O\Bigg(\sum\limits_{k=1}^{n-2}
\big|\psi(k)k^{s}-\psi(k+1)(k+1)^{s}\big|k^{\frac{1}{p}}\Bigg)+
O\left(\psi(n-1)(n-1)^{s+\frac{1}{p}}\right).
\end{equation}

Покажемо, що якщо $g_{s,p}\in A^{+}$, то
\begin{equation}\label{v1a+}
  \Big\|
\sum\limits_{k=1}^{n-1}\psi(k)k^{s}\cos\left(kt+\frac{\beta\pi}{2}\right)\Big\|_{p'}=O\big(\psi(n)n^{s+\frac{1}{p}}\big),
\ 1\leq p<\infty.
\end{equation}
З нерівності трикутника випливає
\begin{equation}\label{v1a+1}
\sum\limits_{k=1}^{n-2}\!
\big|\psi(k)k^{s}-\psi(k+1\!)(k+1)^{s}\big|k^{\frac{1}{p}}\leq\sum\limits_{k=1}^{n-2}\big|\psi(k)k^{s+\frac{1}{p}}-\psi(k+1)(k+1)^{s+\frac{1}{p}}\big|+
$$
$$
+\sum\limits_{k=1}^{n-2}\big|\psi(k+1)(k+1)^{s+\frac{1}{p}}-\psi(k+1)(k+1)^{s}k^{\frac{1}{p}}\big|.
\end{equation}
Далі, використовуючи монотонне зростання функції $g_{s,p}$ з
множини $A^{+}$, отримуємо
\begin{equation}\label{ddzn2}
\sum\limits_{k=1}^{n-2}\!
\big|\psi(k)k^{s}-\psi(k+1\!)(k+1)^{s}\big|k^{\frac{1}{p}}<\sum\limits_{k=1}^{n-2}\big(\psi(k+1)(k+1)^{s+
\frac{1}{p}}-\psi(k)k^{s+\frac{1}{p}}\big)+
$$
$$
+\sum\limits_{k=1}^{n-2}
\psi(k+1)(k+1)^{s}\big((k+1)^{\frac{1}{p}}-k^{\frac{1}{p}}\big)<\psi(n-1)(n-1)^{s+\frac{1}{p}}+
$$
$$
+\frac{1}{p}\sum\limits_{k=1}^{n-2}\psi(k+1)(k+1)^{s}k^{\frac{1}{p}-1}=
$$
$$
=O\Bigg(\psi(n)n^{s+\frac{1}{p}}+
\sum\limits_{k=1}^{n-2}\psi(k+1)(k+1)^{s+\frac{1}{p}-1}\big(1+\frac{1}{k}\big)^{\frac{1}{p'}}\Bigg)=
$$
$$
=O\Bigg(\psi(n)n^{s+\frac{1}{p}}+
\sum\limits_{k=2}^{n-1}\frac{\psi(k)k^{s+\frac{1}{p}}}{k}\Bigg)=
$$
$$
=O\Bigg(\psi(n)n^{s+\frac{1}{p}}+\int\limits_{1}^{n-1}\frac{\psi(t)t^{s+\frac{1}{p}}}{t}dt\Bigg)=O\big(\psi(n)n^{s+\frac{1}{p}}\big).
\end{equation}
Об'єднавши \eqref{ddzn}  і \eqref{ddzn2}, одержуємо \eqref{v1a+}.
  Із \eqref{111}, \eqref{nb} і \eqref{v1a+} при $\psi\in \Theta_{p}$
і $g_{s,p}\in A^{+}$ випливає  \eqref{t5}.

Далі покажемо, що при виконанні умови $g_{s,p}\in A^{-}$ для
величини \eqref{v1} виконується оцінка
\begin{equation}\label{cpadna}
\Bigg\|\sum\limits_{k=1}^{n-1}\psi(k)k^{s}\cos\left(kt+\frac{\beta\pi}{2}\right)\Bigg\|_{p'}=O(1),
\ 1\leq p<\infty.
\end{equation}
В силу \eqref{v1a+1} і монотоного спадання функції $g_{s,p}\in
A^{-}$ отримуємо
\begin{equation}\label{ddzn1}
\sum\limits_{k=1}^{n-2}
\Big|\psi(k)k^{s}-\psi(k+1)(k+1)^{s}\Big|k^{\frac{1}{p}}<
\psi(1)+\frac{1}{p}\sum\limits_{k=1}^{n-2}\psi(k+1)(k+1)^{s}k^{\frac{1}{p}-1}=
$$
$$
=O\Bigg(1+\sum\limits_{k=1}^{n-2}\psi(k+1)(k+1)^{s+\frac{1}{p}-1}\big(1+\frac{1}{k}\big)^{\frac{1}{p'}}\Bigg)=
$$
$$
= O\Bigg(1+
\sum\limits_{k=2}^{n-1}\frac{\psi(k)k^{s+\frac{1}{p}}}{k}\Bigg)=
O\Bigg(1+\int\limits_{1}^{n-1}\frac{\psi(t)t^{s+\frac{1}{p}}}{t}dt\Bigg)=O(1).
\end{equation}
З \eqref{ddzn} та \eqref{ddzn1} одержуємо \eqref{cpadna}.
Очевидно, що умова $g_{s,p}\in A^{-}$ забезпечує включення
$\psi\in \Theta_{p}$. Об'єднавши  \eqref{111}, \eqref{nb} і
\eqref{cpadna} за умови $g_{s,p}\in A^{-}$, отримуємо оцінку
\eqref{t1}.

Для оцінки величини \eqref{v1} у випадку, коли    $g_{s,p}\in
\mathcal{Z}$, нам буде корисним наступне твердження.

\bf Лема 2.   \it{Нехай $0<r<1$ і $g_{s,p}\in \mathcal{Z}$. Тоді
для довільних $x\in (0,\pi]$ і $N\in \mathbb{N}$ мають місце
нерівності
\begin{equation}\label{nsynsun}
\bigg|\sum\limits_{k=1}^{N}\frac{g(k)}{k^{r}}\sin kx\bigg|\leq
C_{g,r}g\Big(\frac{1}{x}\Big)x^{r-1},
\end{equation}
\begin{equation}\label{nsyncos}
\bigg|\sum\limits_{k=1}^{N}\frac{g(k)}{k^{r}}\cos kx\bigg|\leq
\widetilde{C}_{g,r}g\Big(\frac{1}{x}\Big)x^{r-1},
\end{equation}
в яких $C_{g,r}$, $\widetilde{C}_{g,r}$ --- додатні величини, що
залежать лише від $g$ та $r$.

{\textbf{\textit{Доведення.}}} \rm    Покажемо справедливість лише
нерівності \eqref{nsynsun}, оскільки \eqref{nsyncos} доводиться
аналогічно.

Представимо функцію
$S_{g,r}^{(N)}(x)=\sum\limits_{k=1}^{N}\frac{g(k)}{k^{r}}\sin kx$
у вигляді
\begin{equation}\label{predsg}
S_{g,r}^{(N)}(x)=S_{g,r}(x)-\sum\limits_{k=N+1}^{\infty}\frac{g(k)}{k^{r}}\sin
kx,
\end{equation}
де $S_{g,r}(x)$ --- сума ряду
\begin{equation}\label{rydmag}
S_{g,r}(x)=\sum\limits_{k=1}^{\infty}\frac{g(k)}{k^{r}}\sin kx.
\end{equation}
Функція $S_{g,r}$ є неперервною  на $(0,\pi]$, а її частинні суми
$S^{(N)}_{g,r}$ рівномірно збігаються до $S_{g,r}$ на будь-якому
сегменті $[\varepsilon,\pi]$, де $0<\varepsilon<\pi$ (див.,
наприклад, \cite{Z}, c. 6). Крім того, згідно з теоремою (5.2.6)
роботи (\cite{Z}, c. 300), при $x\rightarrow+0$ має місце
асимптотична рівність
\begin{equation}\label{asumriv}
S_{g,r}(x)=g\Big(\frac{1}{x}\Big)x^{r-1}\big(\Gamma(1-r)\cos\frac{\pi
r }{2}
 +o(1)\big), \ 0<r<1,
\end{equation}
в якій $\Gamma(\cdot)$ --- гамма-функція Ейлера. З \eqref{asumriv}
випливає, що
\begin{equation}\label{magor1}
|S_{g,r}(x)|\leq C^{(1)}_{g,r}g\Big(\frac{1}{x}\Big)x^{r-1}, \
0<x\leq\pi,
\end{equation}
 $C^{(1)}_{g,r}$ --- залежить лише від $g$ та
$r$.

Через $\Theta_{g,r}^{(1)}$, $r\in (0,1)$, позначимо  значення
аргументу функції $g(t)t^{-\frac{r}{2}}$ таке, що при $t\in~
[\Theta_{g,r}^{(1)},\infty)$ функція $g(t)t^{-\frac{r}{2}}$
монотонно спадає, а через $\Theta_{g,r}^{(2)}$ --- значення
аргументу функції $g(t)t^{\frac{1-r}{2}}$ таке, що для довільного
$\tau\geq\Theta_{g,r}^{(2)}$
\begin{equation*}
 \max\limits_{t\in
[1,\tau]}g(t)t^{\frac{1-r}{2}}=g(\tau)\tau^{\frac{1-r}{2}},
\end{equation*}
(оскільки $g\in \mathcal{Z}$, то $\Theta_{g,r}^{(i)}$, $i=1,2$
існують для довільного $r\in(0,1)$). Покладемо
$M=M(g,r)=\max\big\{1, \pi\Theta^{(1)}_{g,r},
\pi\Theta^{(2)}_{g,r}\big\}$. Тоді для довільного $x\in (0,\pi]$
\begin{equation}\label{maxt1}
\max\limits_{t\in
[\frac{M}{x},\infty]}g(t)t^{-\frac{r}{2}}=g\Big(\frac{M}{x}\Big)\Big(\frac{x}{M}\Big)^{\frac{r}{2}},
\end{equation}
\begin{equation}\label{maxt2}
\max\limits_{t\in
[1,\frac{M}{x}]}g(t)t^{\frac{1-r}{2}}=g\Big(\frac{M}{x}\Big)\Big(\frac{M}{x}\Big)^{\frac{1-r}{2}}.
\end{equation}

Доведемо справедливість нерівності \eqref{nsynsun} при
$N>\frac{M}{x}$. Застосовуючи перетворення Абеля і враховуючи
\eqref{maxt1}, а також нерівність (5.2.27) роботи (\cite{Z}, с.
306), згідно з якою
\begin{equation*}
\sum\limits_{\nu=1}^{N}\nu^{-\gamma}\sin \nu x\leq
C_{r}x^{\gamma-1}, \ N\in\mathbb{N}, \  \gamma\in (0,1),
\end{equation*}
$C_{r}$ --- залежить лише від $r$, одержимо
\begin{equation}\label{mod1}
\bigg|\sum\limits_{k=N+1}^{\infty}\frac{g(k)\sin kx}{k^{r}}\bigg|=
$$
$$
=\bigg|\sum\limits_{k=N+1}^{\infty}\bigg(\frac{g(k)}{k^{\frac{r}{2}}}-\frac{g(k+1)}{(k+1)^{\frac{r}{2}}}\bigg)
\sum\limits_{\nu=1}^{k}\frac{\sin \nu
x}{\nu^{\frac{r}{2}}}-\frac{g(N+1)}{(N+1)^{\frac{r}{2}}}\sum\limits_{\nu=1}^{N}\frac{\sin
\nu x }{\nu^{\frac{r}{2}}}\bigg|\leq
$$
$$
\leq 2\frac{g(N+1)}{(N+1)^{\frac{r}{2}}}\sup\limits_{k\in
\mathbb{N}}\Big|\sum\limits_{\nu=1}^{k}\frac{\sin \nu x
}{\nu^{\frac{r}{2}}}\Big|\leq 2 C_{r}\max\limits_{N\in
\mathbb{N}}\frac{g(N+1)}{(N+1)^{\frac{r}{2}}}x^{\frac{r}{2}-1}\leq
$$
$$
\leq2C_{r}
g\Big(\frac{M}{x}\Big)\Big(\frac{x}{M}\Big)^{\frac{r}{2}}x^{\frac{r}{2}-1}=C_{g,r}^{(2)}g\Big(\frac{M}{x}\Big)x^{r-1},
\ r\in (0,1].
\end{equation}
Оскільки $g\in \mathcal{Z}$, то (див., наприклад, \cite{Z}, с.
299)
\begin{equation*}
g(Mt)=g(t)+o\big(g(t)\big), \ t\geq1, \ M>0,
\end{equation*}
а, отже,
\begin{equation}\label{slkner}
\widetilde{C}_{M, g}\leq g(Mt)\leq C_{M, g}g(t), \ t\geq1, \ M>0,
\end{equation}
де $\widetilde{C}_{M, g}$, $C_{M,g}$ --- залежать лише від $g$ та
$M$.

В силу \eqref{mod1} та \eqref{slkner} одержуємо
\begin{equation}\label{mod2}
\bigg|\sum\limits_{k=N+1}^{\infty}\frac{g(k)\sin
kx}{k^{r}}\bigg|\leq C_{g,r}^{(3)}g\Big(\frac{1}{x}\Big)x^{r-1}.
\end{equation}
Із \eqref{predsg}, \eqref{magor1} та \eqref{mod2} випливає, що при
$N>\frac{M}{x}$
\begin{equation}\label{nernesk}
\Big|S_{g,r}^{(N)}(x)\Big|\leq
\big|S_{g,r}(x)\big|+\bigg|\sum\limits_{k=N+1}^{\infty}\frac{g(k)}{k^{r}}\sin
kx\bigg|\leq
$$
$$
\leq
C_{g,r}^{(1)}g\Big(\frac{1}{x}\Big)x^{r-1}+C_{g,r}^{(3)}g\Big(\frac{1}{x}\Big)x^{r-1}=C_{g,r}^{(4)}g\Big(\frac{1}{x}\Big)x^{r-1}.
\end{equation}

Переконаємось у справедливості нерівності \eqref{nsynsun} при
$N\leq\frac{M}{x}$. Очевидно, що
\begin{equation}\label{nerm15}
\Big|S_{g,r}^{(N)}(x)\Big|\leq\sum\limits_{k=1}^{N}\frac{g(k)}{k^{r}}=\sum\limits_{k=1}^{N}\frac{g(k)k^{\frac{1-r}{2}}}{k^{\frac{1+r}{2}}}.
\end{equation}
Оскільки при $N\leq\frac{M}{x}$ для $k=\overline{1,N}$ виконується
\eqref{maxt2}, то з урахуванням \eqref{slkner} і \eqref{nerm15}
\begin{equation}\label{nerm2}
\Big|S_{g,r}^{(N)}(x)\Big|\leq
g\Big(\frac{M}{x}\Big)\Big(\frac{M}{x}\Big)^{\frac{1-r}{2}}\sum\limits_{k=1}^{N}\frac{1}{k^{\frac{1+r}{2}}}<
$$
$$
<C_{r}^{(1)}g\Big(\frac{M}{x}\Big)\Big(\frac{M}{x}\Big)^{1-r}=
C_{g,r}^{(5)}g\Big(\frac{1}{x}\Big)x^{r-1}.
\end{equation}

Поєднуючи \eqref{nernesk} та \eqref{nerm2}, одержимо нерівність
\eqref{nsynsun}. Лему 2 доведено.

 За виконання умови $g_{s,p}\in
\mathcal{Z}$, має місце оцінка
 \begin{equation}\label{ozfg}
\sum\limits_{k=1}^{N}\frac{g_{s,p}(k)}{k^{r}}=O\Big(N^{1-r}
g_{s,p}(N)\Big),\ r\in[0,1), \ x\in\mathbb{R},
\end{equation}
яка випливає   з асимптотичної рівності (див., \cite{Z}, с. 299)
\begin{equation}\label{rivdg}
\sum\limits_{k=1}^{N}g_{s,p}(k)k^{\alpha}\sim
\frac{N^{1+\alpha}}{1+\alpha}g_{s,p}(N), \ \alpha>-1,
\end{equation}
при $\alpha=-r$ (запис $A(n)\sim B(n)$ означає виконання
граничного співвідношення
\linebreak$\lim\limits_{n\rightarrow\infty}\frac{A(n)}{B(n)}=1$).

Покажемо, що при $1<p<\infty$ і  $g_{s,p}\in \mathcal{Z}$, має
місце оцінка
 \begin{equation}\label{slabko1}
\Bigg\|\sum\limits_{k=1}^{n-1}\psi(k)k^{s}\cos\left(\!kt+\frac{\beta\pi}{2}\!\right)\!\Bigg\|_{p'}\!\!=
O\Bigg(\int\limits_{1}^{n}\frac{(\psi(\tau)\tau^{s+\frac{1}{p}})^{p'}}{\tau}d\tau\Bigg)^{\frac{1}{p'}}.
\end{equation}

 Враховуючи нерівності
\eqref{nsynsun}, \eqref{nsyncos} та \eqref{ozfg}, при $1<p<\infty$
одержимо
\begin{equation}\label{slabko12}
\Bigg\|\sum\limits_{k=1}^{n-1}\psi(k)k^{s}\cos\left(\!kt+\frac{\beta\pi}{2}\!\right)\!\Bigg\|_{p'}=
$$
$$
=O\left(\Bigg(\int\limits_{0\leq\tau\leq\frac{1}{n}}
\Bigg(\sum\limits_{k=1}^{n-1}\frac{\psi(k)k^{s+\frac{1}{p}}}{k^{\frac{1}{p}}}\Bigg)^{p'}\!\!d\tau+\int\limits_{\frac{1}{n}\leq\tau\leq\pi}\Bigg(
\frac{\psi\big(\frac{1}{\tau}\big)\tau^{\frac{1}{p}-1}}{\tau^{s+\frac{1}{p}}}\Bigg)^{p'}d\tau\Bigg)^{\frac{1}{p'}}\right)=
$$
$$
=
O\left(\Bigg(\Big(\psi(n)n^{s+\frac{1}{p}}\Big)^{p'}\big(n^{1-\frac{1}{p}}\big)^{p'}n^{-1}+\int\limits_{\frac{1}{\pi}\leq\tau\leq
n }\Big(\frac{\psi(
\tau)\tau^{s+\frac{1}{p}}}{\tau^{\frac{1}{p}-1}}\Big)^{p'}\frac{d\tau}{\tau^{2}}\Bigg)^{\frac{1}{p'}}\right)=
$$
$$
=O\Bigg(\Big(\psi(n)n^{s+\frac{1}{p}}\Big)^{p'}+\int\limits_{\frac{1}{\pi}\leq\tau\leq
n}
\Big(\frac{\psi(\tau)\tau^{s+\frac{1}{p}}}{\tau^{-\frac{1}{p'}}}\cdot\frac{1}{\tau^{\frac{2}{p'}}}\Big)^{p'}d\tau\Bigg)^{\frac{1}{p'}}=
$$
$$
=O\Bigg(\Big(\psi(n)n^{s+\frac{1}{p}}\Big)^{p'}+\int\limits_{1}^{n}\frac{(\psi(\tau)\tau^{s+\frac{1}{p}})^{p'}}
{\tau}d\tau\Bigg)^{\frac{1}{p'}}=
$$
$$
=O\Bigg(g^{p'}_{s,p}(n)+\int\limits_{1}^{n}\frac{g^{p'}_{s,p}(\tau)}{\tau}d\tau\Bigg)^{\frac{1}{p'}}.
\end{equation}

Зауважимо, що якщо $g_{s,p}\in \mathcal{Z}$, то  $g_{s,p}^{p'}\in
\mathcal{Z}$. Оскільки $g^{p'}_{s,p}\in \mathcal{Z}$, то має місце
оцінка  (див., наприклад, \cite{Z}, с. 302)
\begin{equation}\label{normap''}
g_{s,p}^{p'}(n)=O\Big(\int\limits_{1}^{n}\frac{g_{s,p}^{p'}(t)}{t}dt\Big),
\ n\in\mathbb{N}.
\end{equation}
Із  \eqref{slabko12} і \eqref{normap''}    випливає
\eqref{slabko1}.

Неважко переконатись, що умова $g_{s,p}\in\mathcal{Z}$ при $s>0$ і
$1<p<\infty$ забезпечує включення $\psi\in \Theta_{p}$. Тому з
\eqref{111}, \eqref{nb} та \eqref{slabko1} при $1<p<\infty$ і
$g_{s,p}\in\mathcal{Z}$ отримуємо \eqref{t4}.
 Теорему 1 доведено.

\bf Теорема 2.  \it{Нехай $1< p<\infty$, $s>0$,
$g_{s,p}(t):=\psi(t)t^{s+\frac{1}{p}}$. Тоді

1. Якщо $\psi\in \Theta_{p}$ і  $g_{s,p}\in A^{+}$,  то для
довільних
 $\beta\in \mathbb{R}$ і $n\in \mathbb{N}$
\begin{equation}\label{t52}
{\cal E}_{n}\left(C^{\psi}_{\beta,p};
Z_{n}^{s}\right)_{C}\asymp{E}_{n}(C^{\psi}_{\beta,p})_{C}\asymp\psi(n)n^{\frac{1}{p}}.
\end{equation}

2. Якщо   $g_{s,p}\in \mathcal{Z}^{+}_{p'}$ або $g_{s,p}\in
\mathcal{Z}^{-}_{p'}$ і
$\int\limits_{1}^{n}\frac{g^{p'}_{s,p}(t)}{t}dt\neq O(1)$, то для
довільних $\beta \in \mathbb{R}$ і $n\in \mathbb{N}$
\begin{equation}\label{t42}
{\cal E}_{n}\left(C^{\psi}_{\beta,p};
Z_{n}^{s}\right)_{C}\asymp\frac{1}{n^{s}}\Bigg(\int\limits_{1}^{n}\frac{\big(\psi(t)t^{s+\frac{1}{p}}\big)^{p'}}{t}dt\Bigg)
^{\frac{1}{p'}}, \ \frac{1}{p}+\frac{1}{p'}=1.
\end{equation}

3. Якщо $g_{s,p}\in A^{-}$ або  $g_{s,p}\in \mathcal{Z}^{-}_{p'}$
і $\int\limits_{1}^{n}\frac{g^{p'}_{s,p}(t)}{t}dt=O(1)$,
$\frac{1}{p}+\frac{1}{p'}=1$, то для довільних $\beta \in
\mathbb{R}$ і $n\in \mathbb{N}$
\begin{equation}\label{t12}
{\cal E}_{n}\left(C^{\psi}_{\beta,p};
Z_{n}^{s}\right)_{C}\asymp\frac{1}{n^{s}}.
\end{equation}

{\textbf{\textit{Доведення.}}} \ \rm Оцінки зверху у \eqref{t52}
-- \eqref{t12} випливають із співвідношень \eqref{t5} --
\eqref{t1} теореми 1. Встановимо необхідні оцінки знизу величин
${\cal E}_{n}\left(C^{\psi}_{\beta,p};Z_{n}^{s}\right)_{C}$ .

Позначимо через $B$ множину монотонно незростаючих додатних
функцій $\psi(t)$, заданих на $[1,\infty)$, для кожної з яких
існує стала $K>0$, така, що
\begin{equation}\label{umb}
\frac{\psi(t)}{\psi(2t)}\leq K, \ t\geq 1.
\end{equation}
У роботі (\cite{Gra}, с. 4) при $1<p<\infty$ і $\psi\in
B\cap\Theta_{p}$ встановлено  оцінку знизу для
${E}_{n}(C^{\psi}_{\beta,p})_{C}$:
\begin{equation}\label{poroznnn}
{E}_{n}(C^{\psi}_{\beta,p})_{C}\geq
K_{\psi,p}\psi(n)n^{\frac{1}{p}}, \ n\in \mathbb{N}, \ \beta\in
\mathbb{R},
\end{equation}
де $K_{\psi,p}$ --- додатні величини, що можуть залежати лише від
$\psi$ та $p$. Оскільки
\begin{equation}\label{nerznabl}
{E}_{n}(C^{\psi}_{\beta,p})_{C}\leq{\cal
E}_{n}\left(C^{\psi}_{\beta,p};Z_{n}^{s}\right)_{C}, \ n\in
\mathbb{N},
\end{equation}
то з нерівностей \eqref{poroznnn}, \eqref{nerznabl} та оцінки
\eqref{t5}  випливає порядкова рівність \eqref{t52}, якщо буде
доведена імплікація
\begin{equation}\label{vrl1}
g_{s,p}\in A^{+}\Rightarrow \psi\in B, \ s>0, \ 1< p<\infty.
\end{equation}

Якшо $g_{s,p}\in A^{+}$, то існує $\varepsilon>0$ і монотонно
неспадна на $[1,\infty)$ функція $\varphi$ така, що
$\psi(t)=\frac{\varphi(t)}{t^{s+\frac{1}{p}-\varepsilon}}$, звідки
\begin{equation*}
\frac{\psi(t)}{\psi(2t)}=\frac{\varphi(t)}{\varphi(2t)}2^{s+\frac{1}{p}-\varepsilon}\leq
2^{s+\frac{1}{p}-\varepsilon}=K,
\end{equation*}
і, отже, $\psi\in B$.

Втім, оцінку знизу величини ${\cal E}_{n}\left(C^{\psi}_{\beta,p};
Z_{n}^{s}\right)_{C}$ при $1<p<\infty$, \linebreak$\psi\in
\Theta_{p}$ та $g_{s,p}\in A^{+}$ можна одержати також, якщо
розглянути функцію $f_{1}(t)=(\Psi_{\beta}\ast\varphi_{1})(t)$, де
\begin{equation*}
\varphi_{1}(t)=\frac{a_{1}}{n^{\frac{1}{p'}}}\sum\limits_{k=1}^{n-1}\cos\big(kt+\frac{\beta\pi}{2}\big),
\ \frac{1}{p}+\frac{1}{p'}=1,
\end{equation*}
 $a_{1}$ --- додатний параметр, значення  якого буде вказано
пізніше. В силу \eqref{ndiri} маємо
\begin{equation}\label{fi1}
\left\Vert\varphi_{1}\right\Vert_{p}=\frac{a_{1}}{n^{\frac{1}{p'}}}\left\Vert\sum\limits_{k=1}^{n-1}\cos\big(kt+\frac{\beta\pi}{2}\big)\right\Vert_{p}\leq
\frac{a_{1}K_{1}n^{\frac{1}{p'}}}{n^{\frac{1}{p'}}}=a_{1}K_{1}, \
1<p<\infty.
\end{equation}
При $a_{1}=(K_{1})^{-1}$ з \eqref{fi1} одержимо  нерівність
$\left\Vert\varphi_{1}\right\Vert_{p}\leq1$, а, отже,  включення
$f_{1}\in C^{\psi}_{\beta,p}$.

Як випливає з формули (1.9) з монографії (\cite{Z}, с. 65), для
$f_{1}$ має місце рівність
\begin{equation*}
f_{1}(t)=\frac{a_{1}}{n^{\frac{1}{p'}}}\sum\limits_{k=1}^{n-1}\psi(k)\cos
kt.
\end{equation*}
Тоді для функції $f_{1}(t)$ справедлива оцінка
\begin{equation}\label{ochzndonas}
{\cal E}_{n}\left(C^{\psi}_{\beta,p};
Z_{n}^{s}\right)_{C}\geq\big|f_{1}(0)-Z_{n}^{s}(f_{1},0)\big|=\frac{a_{1}}{n^{s+\frac{1}{p'}}}\sum\limits_{k=1}^{n-1}\psi(k)k^{s}\geq
$$
$$
\geq\frac{a_{2}\psi(n)}{n^{s+\frac{1}{p'}}}\sum\limits_{k=1}^{n-1}k^{s}\geq
K_{2}\psi(n)n^{\frac{1}p}.
\end{equation}

Доведемо оцінку знизу величини ${\cal
E}_{n}\left(C^{\psi}_{\beta,p}; Z_{n}^{s}\right)_{C}$ при
$1<p<\infty$ і   $g_{s,p}\in \mathcal{Z}^{+}_{p'}$. Для цього
розглянемо функцію
\begin{equation}\label{f}
\varphi_{2}(\tau)\!=\!\!\frac{a_{2}}{\Bigg(\sum\limits_{k=1}^{n-1}\frac{\big(\psi(k)k^{s+\frac{1}{p}}\big)^{p'}}{k}\Bigg)^{\frac{1}{p}}}\sum\limits_{k=1}^{n-1}
\frac{\big(\psi(k)k^{s+\frac{1}{p}}\big)^{p'-1}}{k^{\frac{1}{p'}}}\cos\Big(\!k\tau+\frac{\beta\pi}{2}\Big),
\ a_{2}>0.
\end{equation}
 Покладемо
\begin{equation*}
 f_{2}(\tau)=(\Psi_{\beta}\ast\varphi_{2})(\tau)=
$$
$$
=\frac{a_{2}}
{\Bigg(\sum\limits_{k=1}^{n-1}\frac{\big(\psi(k)k^{s+\frac{1}{p}}\big)^{p'}}{k}\Bigg)^{\frac{1}{p}}}\Bigg
(\sum\limits_{k=1}^{n-1}\frac{\psi^{p'}(k)(k^{s+\frac{1}{p}})^{p'-1}}{k^{\frac{1}{p'}}}\cos
k\tau\Bigg),
\end{equation*}
і покажемо, що при відповідному  підборі параметра  $a_{2}$
виконуватиметься включення \linebreak$f_{2}\in
C^{\psi}_{\beta,p}$.

 Застосувавши  перетворення Абеля до суми
\begin{equation*}
 \sum\limits_{k=1}^{n-1}
\frac{\big(\psi(k)k^{s+\frac{1}{p}}\big)^{p'-1}}{k^{\frac{1}{p'}}}\cos\Big(k\tau+\frac{\beta\pi}{2}\Big),
\end{equation*}
одержимо
\begin{equation}\label{n1}
\Bigg\|\sum\limits_{k=1}^{n-1}
\frac{\big(\psi(k)k^{s+\frac{1}{p}}\big)^{p'-1}}{k^{\frac{1}{p'}}}\cos\Big(k\tau+\frac{\beta\pi}{2}\Big)\Bigg\|_{p}=
$$
$$
=\!\Bigg\|\sum\limits_{k=1}^{n-2}\Big(\big(\psi(k)k^{s+\frac{1}{p}}\big)^{p'-1}\!-\!
\big(\psi(k+1)(k+1)^{s+\frac{1}{p}}\big)^{p'-1}\Big)\!\!\sum\limits_{\nu=1}^{k}\frac{\cos\big(\nu\tau+\frac{\beta\pi}{2}\big)}
{\nu^{\frac{1}{p'}}}+
$$
$$
+\big(\psi(n-1)(n-1)^{s+\frac{1}{p}}\big)^{p'-1}\sum\limits_{\nu=1}^{n-1}\frac{\cos\big(\nu\tau+\frac{\beta\pi}{2}\big)}{\nu^{\frac{1}{p'}}}\Bigg\|_{p}
\leq
$$
$$
\leq\sum\limits_{k=1}^{n-2}\big|\big(\psi(k)k^{s+\frac{1}{p}}\big)^{p'-1}\!-\!
\big(\psi(k+1)(k+1)^{s+\frac{1}{p}}\big)^{p'-1}\big|\Big\|\sum\limits_{\nu=1}^{k}\frac{\cos\big(\nu\tau+\frac{\beta\pi}{2}\big)}{\nu^{\frac{1}{p'}}}\Big\|_{p}+
$$
$$
+\big(\psi(n-1)(n-1)^{s+\frac{1}{p}}\big)^{p'-1}\Big\|\sum\limits_{\nu=1}^{n-1}\frac{\cos\big(\nu\tau+\frac{\beta\pi}{2}\big)}{\nu^{\frac{1}{p'}}}\Big\|_{p}.
\end{equation}

Оскільки (див., наприклад, \cite{Z}, с. 306)
\begin{equation*}
\Bigg|\sum\limits_{k=1}^{n}\frac{\cos
k\tau}{k^{\frac{1}{p}}}\Bigg|={\left\{ {\begin{array}{l l}
 O\big(n^{1-\frac{1}{p}}\big), & p>1, \ \tau\in\mathbb{R}, \\
O\big(|\tau|^{\frac{1}{p}-1}\big), & 0<|\tau|\leq\pi, \ p>1,
\end{array}} \right.}
$$
$$
\Bigg|\sum\limits_{k=1}^{n}\frac{\sin
k\tau}{k^{\frac{1}{p}}}\Bigg|={\left\{ {\begin{array}{l l}
 O\big(n^{1-\frac{1}{p}}\big), & p>1, \ \tau\in\mathbb{R}, \\
O\big(|\tau|^{\frac{1}{p}-1}\big), & 0<|\tau|\leq\pi, \ p>1,
\end{array}} \right.}
\end{equation*}
то   при $1< p<\infty$,  $n\in\mathbb{N}$ і $\beta\in\mathbb{R}$
\begin{equation}\label{nk}
\Bigg\|\sum\limits_{k=1}^{n}\frac{\cos\big(kt+\frac{\beta\pi}{2}\big)}{k^{\frac{1}{p'}}}\Bigg\|_{p}\leq\widetilde{K}_{p}^{(1)}\ln^{\frac{1}{p}}n,
\end{equation}
де $\widetilde{K}_{p}^{(1)}$ --- стала, що залежать лише від $p$.
З урахуванням  оцінок \eqref{f} --- \eqref{nk}  будемо мати
\begin{equation}\label{nfi}
\|\varphi_{2}(\tau)\|_{p}\leq\frac{a_{2}\widetilde{K}_{p}^{(1)}}{\Bigg(\sum\limits_{k=1}^{n-1}\frac{\big(\psi(k)k^{s+\frac{1}{p}}\big)^{p'}
}{k}\Bigg)^{\frac{1}{p}}} \times
$$
$$
\times\Big(\sum\limits_{k=1}^{n-2}|(\psi(k)k^{s+\frac{1}{p}})^{p'-1}-
(\psi(k+1)(k+1)^{s+\frac{1}{p}})^{p'-1}|\ln^{\frac{1}{p}}k+
$$
$$
+(\psi(n-1)(n-1)^{s+\frac{1}{p}})^{p'-1}\ln^{\frac{1}{p}}(n-1)\Big).
\end{equation}
 В силу
монотонного неспадання функції $g_{s,p}$ із множини
$\mathcal{Z}^{+}_{p'}$, одержимо
\begin{equation}\label{ch}
\sum\limits_{k=1}^{n-2}\big|\big(\psi(k)k^{s+\frac{1}{p}}\big)^{p'-1}\!-\!
\big(\psi(k+1)(k+1)^{s+\frac{1}{p}}\big)^{p'-1}\big|\ln^{\frac{1}{p}}k\leq
$$
$$
\leq \big(\psi(n)n^{s+\frac{1}{p}}\big)^{p'-1}\ln^{\frac{1}{p}}n.
\end{equation}
Із \eqref{nfi} і \eqref{ch} випливає, що за умови $g_{s,p}\in
\mathcal{Z}^{+}_{p'}$
\begin{equation}\label{ochg}
\|\varphi_{2}(\tau)\|_{p}\leq\!\frac{a_{2}\widetilde{K}_{p}^{(1)}g_{s,p}^{p'-1}(n)
\ln^{\frac{1}{p}}n}{\bigg(\sum\limits_{k=1}^{n-1}\frac{g^{p'}(k)}{k}\bigg)^{\frac{1}{p}}}\leq\!
\Bigg(\frac{a_{2}^{p}\widetilde{K}_{p}^{(2)}g^{p'}(n) \ln n
}{\int\limits_{1}^{n}\frac{g^{p'}(t)}{t}dt}\Bigg)^{\frac{1}{p}}\!\!\leq\!\!\big(a_{2}^{p}\widetilde{K}_{p}^{(3)}\big)^{\frac{1}{p}}.
\end{equation}
Вибравши параметр
$a_{2}=(\widetilde{K}_{p}^{(3)})^{-\frac{1}{p}}$, одержимо
нерівність $\|\varphi_{2}\|_{p}\leq 1$, а разом із нею включення
$f_{2}\in C^{\psi}_{\beta,p}$. Для функції $f_{2}$  має місце
оцінка
\begin{equation}\label{to}
{\cal E}_{n}\left(C^{\psi}_{\beta,p};
Z_{n}^{s}\right)_{C}\geq|f_{2}(0)-Z_{n}^{s}(f_{2};0)|=
$$
$$
=
\frac{a_{2}}{n^{s}\Bigg(\sum\limits_{k=1}^{n-1}\frac{\big(\psi(k)k^{s+\frac{1}{p}}\big)^{p'}}{k}\Bigg)^{\frac{1}{p}}}
\sum\limits_{k=1}^{n-1}\frac{\big(\psi(k)
k^{s+\frac{1}{p}}\big)^{p'}}{k}=
$$
$$
=\frac{a_{2}}{n^{s}}\Bigg(\sum\limits_{k=1}^{n-1}\frac{\big(\psi(k)k^{s+\frac{1}{p}}\big)^{p'}}{k}\Bigg)^{\frac{1}{p'}}
\geq\frac{\widetilde{K}_{p}^{(4)}}{n^{s}}\Bigg(\int\limits_{1}^{n}\frac{\big(\psi(\tau)\tau^{s+\frac{1}{p}}\big)^{p'}}
{\tau}d\tau\Bigg)^{\frac{1}{p'}}.
\end{equation}
В силу \eqref{t4} і  \eqref{to} порядкова оцінка \eqref{t42}
доведена за умови $g_{s,p}\in\mathcal{Z}^{+}_{p'}$.

Покажемо, що ця ж оцінка має місце і у випадку
$g_{s,p}\in\mathcal{Z}^{-}_{p'}$ і
$\int\limits_{1}^{n}\frac{g^{p'}_{s,p}(t)}{t}dt\neq O(1)$.
Розглянемо функцію
\begin{equation*}
\varphi_{3}(t)=\frac{a_{3}}{\ln^{\frac{1}{p}}n}\sum\limits_{k=1}^{n-1}\frac{\cos\big(kt+\frac{\beta\pi}{2}\big)}{k^{\frac{1}{p'}}},
\  a_{3}>0,
\end{equation*}
і покладемо
\begin{equation*}
f_{3}(t)=(\varphi_{3}\ast\Psi_{\beta})(t)=\frac{a_{3}}{\ln^{\frac{1}{p}}n}\sum\limits_{k=1}^{n-1}\frac{\psi(k)}{k^{\frac{1}{p'}}}\cos
kt.
\end{equation*}
Покажемо, що параметр $a_{3}$ можна підібрати так, щоб $f_{3}\in
C^{\psi}_{\beta,p}$. Дійсно, в силу \eqref{nk}   при $1<p<\infty$
\begin{equation*}
\|\varphi_{3}\|_{p}=\frac{a_{3}}{\ln^{\frac{1}{p}}n}\Bigg\|\sum\limits_{k=1}^{n-1}\frac{\cos\big(kt+\frac{\beta\pi}{2}\big)}{k^{\frac{1}{p'}}}\Bigg\|_{p}
\leq
\frac{a_{3}\widetilde{K}^{(1)}_{p}\ln^{\frac{1}{p}}n}{\ln^{\frac{1}{p}}n}=a_{3}\widetilde{K}^{(1)}_{p},
\end{equation*}
і тому при     $a_{3}=(\widetilde{K}^{(1)}_{p})^{-1}$, одержимо
нерівність
 $\|\varphi_{3}\|_{p}\leq 1$, а разом з нею і  включення  $f_{3}\in C^{\psi}_{\beta,p}$.

Для $f_{3}$ за умови  $g_{s,p}\in \mathcal{Z}^{-}_{p'}$ має місце
оцінка
\begin{equation}\label{ff6}
{\cal E}_{n}\left(C^{\psi}_{\beta,p};
Z_{n}^{s}\right)_{C}\geq|f_{3}(0)-Z^{s}_{n}(f_{3};0)|=\frac{a_{3}}{n^{s}\ln^{\frac{1}{p}}n}\sum\limits_{k=1}^{n-1}\frac{\psi(k)k^{s+\frac{1}{p}}}{k}=
$$
$$
=\frac{a_{3}}{n^{s}\ln^{\frac{1}{p}}n}\sum\limits_{k=1}^{n-1}\frac{g_{s,p}(k)}{k}\geq
\frac{a_{3}g_{s,p}(n)}{n^{s}\ln^{\frac{1}{p}}n}
\sum\limits_{k=1}^{n-1}\frac{1}{k}\geq
\frac{\widetilde{K}_{p}^{(5)}g_{s,p}(n)\ln^{\frac{1}{p'}}n}{n^{s}}=
$$
$$
=\frac{\widetilde{K}_{p}^{(5)}\big(g^{p'}_{s,p}(n)\ln n
\big)^{\frac{1}{p'}}}{n^{s}}
\geq\frac{\widetilde{K}_{p}^{(6)}}{n^{s}}\int\limits_{1}^{n}\Big(\frac{g^{p'}_{s,p}(t)}{t}dt\Big)^{\frac{1}{p'}}.
\end{equation}
Із  \eqref{t4} та \eqref{ff6} одержуємо порядкову рівність
\eqref{t42} при $1<p<\infty$ та $g_{s,p}\in \mathcal{Z}^{-}_{p'}$
і \linebreak$\int\limits_{1}^{n}\frac{g^{p'}_{s,p}(t)}{t}dt\neq
O(1)$.

 Як випливає  з теореми (2.2.1) роботи
\cite{S1}, с. 92,  метод $Z_{n}^{s}$ насичений з порядком
насичення ${n^{-s}}$, а це означає, що величина ${\cal
E}_{n}\left(C^{\psi}_{\beta,p}; Z_{n}^{s}\right)_{C}$ не може
спадати до нуля  швидше, ніж ${n^{-s}}$, тобто має місце оцінка
\begin{equation}\label{monnas}
{\cal E}_{n}\left(C^{\psi}_{\beta,p};
Z_{n}^{s}\right)_{C}\geq\frac{K_{3}}{n^{s}}.
\end{equation}

 Утім, в справедливість нерівності \eqref{monnas} можна преконатись безпосередньо, розглянувши функцію
\begin{equation*}
f_{4}(t)=\frac{\psi(1)}{(2\pi)^{\frac{1}{p}}}\cos t.
\end{equation*}
Легко бачити, що для довільного $\beta \in \mathbb{R}$ і $1\leq
p<\infty$, $f_{4}\in C^{\psi}_{\beta,p}$ і
\begin{equation*}
Z^{s}_{n}(f_{4},t)=\big(1-\frac{1}{n^{s}}\big)\Big(\frac{\psi(1)}{(2\pi)^{\frac{1}{p}}}\cos
t\Big).
\end{equation*}
 Тому
\begin{equation}\label{mon02}
{\cal E}_{n}\left(C^{\psi}_{\beta,p};
Z_{n}^{s}\right)_{C}\geq\big|f_{4}(0)-Z^{s}_{n}(f_{4},0)\big|=\frac{\psi(1)}{(2\pi)^{\frac{1}{p}}
n^{s}}=\frac{K_{3}}{n^{s}}.
\end{equation}
З нерівності \eqref{monnas} та оцінок \eqref{t4}, \eqref{t1}
випливає, що для довільного $\beta\in\mathbb{R}$, $1< p<\infty$ за
умови $g_{s,p}\in \mathcal{Z}^{-}_{p'}$ і
$\int\limits_{1}^{n}\frac{g^{p'}_{s,p}(t)}{t}dt=O(1)$ або
$g_{s,p}\in A^{-}$ виконується порядкова рівність \eqref{t12}.

 Теорему 2
доведено.

Зауважимо, що при  $1< p<\infty$, $s>0$ у випадку, коли $g_{s,p}$
--- рівномірно обмежена зверху і знизу монотонна  функція, як випливає з \eqref{t42},
  для довільних
$\beta\in \mathbb{R}$ і $n\in \mathbb{N}$ виконується порядкова
рівність
\begin{equation}\label{nasl3}
{\cal E}_{n}\left(C^{\psi}_{\beta,p};
Z_{n}^{s}\right)_{C}\asymp\frac{\ln^{1-\frac{1}{p}}n}{n^{s}}.
\end{equation}

Як зазначалось вище, при $s=1$ суми Зигмунда $Z_{n}^{s}$
співпадають із сумами Фейєра $\sigma_{n}$. Тому із теореми 2
випливає наступне  твердження.

\bf Наслідок 1.  \it{Нехай $1< p<\infty$,
 $g_{p}(t):=\psi(t)t^{1+\frac{1}{p}}$. Тоді

1. Якщо $\psi\in \Theta_{p}$ і  $g_{p}\in A^{+}$,  то для
довільних  $\beta\in \mathbb{R}$ і $n\in \mathbb{N}$
\begin{equation}\label{t52f2}
{\cal E}_{n}\left(C^{\psi}_{\beta,p};
\sigma_{n}\right)_{C}\asymp{E}_{n}(C^{\psi}_{\beta,p})_{C}\asymp\psi(n)n^{\frac{1}{p}}.
\end{equation}

2. Якщо  $g_{p}\in \mathcal{Z}^{+}_{p'}$ або $g_{p}\in
\mathcal{Z}^{-}_{p'}$ і
$\int\limits_{1}^{n}\frac{g^{p'}_{p}(t)}{t}dt\neq O(1)$, то для
довільних $\beta \in \mathbb{R}$ і $n\in \mathbb{N}$
\begin{equation}\label{t42f}
{\cal E}_{n}\left(C^{\psi}_{\beta,p};
\sigma_{n}\right)_{C}\asymp\frac{1}{n}\Bigg(\int\limits_{1}^{n}\frac{\big(\psi(t)t^{1+\frac{1}{p}}\big)^{p'}}{t}dt\Bigg)^{\frac{1}{p'}},
\ \frac{1}{p}+\frac{1}{p'}=1.
\end{equation}

3. Якщо $g_{p}\in A^{-}$ або  $g_{p}\in \mathcal{Z}^{-}_{p'}$ і
$\int\limits_{1}^{n}\frac{g^{p'}_{p}(t)}{t}dt=O(1)$,
$\frac{1}{p}+\frac{1}{p'}=1$, то для довільних $\beta \in
\mathbb{R}$ і $n\in \mathbb{N}$
\begin{equation}\label{t12f}
{\cal E}_{n}\left(C^{\psi}_{\beta,p};
\sigma_{n}\right)_{C}\asymp\frac{1}{n}.
\end{equation}
\rm При $\psi(k)=k^{-r}$, $\beta\in \mathbb{R}$ і $r>\frac{1}{p}$
з рівностей \eqref{t52}, \eqref{t12} та \eqref{nasl3} одержуємо
твердження.

\bf Наслідок 2. \it{Нехай $1<p<\infty$, $\beta\in \mathbb{R}$,
$s>0$ і   $n\in
 \mathbb{N}$. Тоді
\begin{equation}\label{vnad}
{\cal E}_{n}\left(W^{r}_{\beta,p};
Z_{n}^{s}\right)_{C}\asymp{\left\{ {\begin{array}{l l}
n^{-r+\frac{1}{p}}, & \frac{1}{p}<r<s+\frac{1}{p}; \\
\frac{\ln^{1-\frac{1}{p}}n}{n^{s}}, &r=s+\frac{1}{p}; \\
n^{-s}, &r>s+\frac{1}{p}.
\end{array}} \right.}
\end{equation}
\rm Оцінки \eqref{vnad} було отримано в роботі
   \cite{Kostuch1}.



\renewcommand{\refname}{}
\makeatletter\renewcommand{\@biblabel}[1]{#1.}\makeatother

\end{document}